\documentclass[11pt,draft]{article}

\font\azhur=msbm10
\font\got=eufm10

\def\smdr{\mbox{\azhur\char110}}

\def\R{\mathop{\hbox{\azhur R}}\nolimits}
\def\N{\mathop{\hbox{\azhur N}}\nolimits}
\def\al{\alpha}
\def\be{\beta}
\def\de{\delta}
\def\ep{\varepsilon}
\def\vf{\varphi}
\def\ga{\gamma}
\def\la{\lambda}
\def\La{\Lambda}

\def\ka{\kappa}

\def\Cup{\mathop{\mathchar4947}\limits}
\def\iff{if and only if}
\def\ccm{Carnot--Caratheo\-dory metric}

\def\Int{\mathop{\rm Int}\nolimits}
\def\Ad{\mathop{\rm Ad}\nolimits}
\def\ad{\mathop{\rm ad}\nolimits}
\def\cont{\mathop{\rm cont}\nolimits}

\def\dist{\mathop{\rm dist}\nolimits}
\def\gT{{\mbox{\got T}}}
\def\gR{{\mbox{\got R}}}
\def\cB{{\cal B}}
\def\cE{{\cal E}}
\def\cF{{\cal F}}
\def\cG{{\cal G}}
\def\cL{{\cal L}}
\def\cN{{\cal N}}
\def\cZ{{\cal Z}}
\def\ins#1{\quad\mbox{#1}\quad}
\def\inv#1{{#1}^{-1}}
\def\iff{if and only if }
\def\empty{\mathord{\hbox{\azhur ?}}}
\def\abs#1{\vert{#1}\vert}
\def\norm#1{\Vert{#1}\Vert}
\def\fracm#1#2{{\textstyle{{#1}\over{#2}}}}
\let\sssize=\scriptscriptstyle
\def\mult#1{\mathop{\cdot_{\sssize #1}}\nolimits}
\let\td=\tilde
\let\wt=\widetilde
\let\ov=\overline
\def\tskl{\vskip0.3\baselineskip}
\def\proof{\vskip0.3\baselineskip\it Proof.\kern10pt\rm}
\def\qed{\kern7pt{$\mathsurround=0pt\bullet$}\par\tskl\rm}
\def\qedm{{.\hfill\kern7pt\bullet}}

\newcounter{xxx}[section]
\newcommand{\lasteqn}[1]%
{\setcounter{xxx}{\value{equation}}\addtocounter{xxx}{-#1}
(\arabic{xxx})}
\newcommand{\lasttheo}[1]%
{\setcounter{xxx}{\value{theorem}}\addtocounter{xxx}{-#1}
\arabic{xxx}}
\newcommand{\lastlem}[1]%
{\setcounter{xxx}{\value{lemma}}\addtocounter{xxx}{-#1}
\arabic{xxx}}

\newcounter{ss}[section]
\newcommand{\ssun}[1]{{\par\noindent\bf\addtocounter{ss}{1}
\arabic{section}.\arabic{ss}~~#1.\kern5pt}}
\newtheorem{lemma}{Lemma}
\newtheorem{theorem}{Theorem}
\newtheorem{corollary}{Corollary}
\title{Attainable  sets for left invariant control systems
and Carnot--Caratheodory metrics on nilpotent Lie groups}
\renewcommand{\thefootnote}{\fnsymbol{footnote}}
\author{V.M.\,Gichev
\thanks{The author was supported in part by INTAS Grant No. 97/10170 and
Federal programm "Integratsia" grant No. 586.}
}
\begin{document}
%
%
\renewcommand{\thefootnote}{\arabic{footnote}}

\date{}
\maketitle

\begin{abstract}
Let $G$ be a semidirect product of a nilpotent Lie group and $\R$.
For a left invariant control system on $G$ whose
control domain is defined by a convex cone in the Lie algebra of $G$ we
prove that the attainable  set coincides with a "halfspace" if the degree
of contact of the cone with certain linear subspaces of the Lie algebra
is sufficiently high.
\end{abstract}

\section{Introduction}
Let $N$ be a Lie group and $C$ be a subset of its Lie algebra $\cN$
identified with the tangent space to $N$ at the identity $e$.
Denote by $\gT(C)$ the set of all piecewise smooth curves with both one-side
tangent vectors in the corresponding left translation of $C$:
\begin{equation}\label{tan}
\ga'(t)\in d_e\la_{\ga(t)}(C)\ins{where}\la_g(h)=gh.
\end{equation}
The attainable  set $\gR(C,p)$ is the closure
of endpoints for curves in $\gT(C)$ which start at $p$.
Put $\gR(C)=\gR(C,e)$ where $e$ is the identity. We consider the problem:
given $C$, find $\gR(C)$. In this setting we may assume without loss of
generality that $C$ is a closed convex cone. If $\gR(C)=N$ then $C$
is called {\it controllable}. The opposite property is the globality:
put
$$H(C)=\{\xi\in\cN:\,\exp(t\xi)\in\gR(C)\ins{for all}t\geq0\},$$
then $C$ is called {\it global} if $H(C)=C$. 
The group $N$ is supposed to be nilpotent and simply connected.
For generating cones in nilpotent Lie algebras the criterion
of controllability due to Hilgert,
Hofmann, and Lawson \cite{HHL}, is known from early 80-th:
\begin{equation}\label{ccrit}
\gR(C)=N\ins{\iff}\Int C\cap[\cN,\cN]\neq\empty.
\end{equation}
If $C\cap[\cN,\cN]=\{0\}$ then the cone is global.
We consider the intermediate case
$$\Int(C)\cap[\cN,\cN]=\empty,\ins{}C\cap[\cN,\cN]\neq\{0\}.$$
This situation can be described as follows: there exist a boundary point
of $C$ and a supporting hyperplane $H$ at this point which
includes  $[\cN,\cN]$. Then the set $\gR(C)$ cannot coincide
with $N$ -- it is included to a "halfspace"
$$N^+_\chi=\inv\chi(\R^+),$$
where $\chi$ is the continuous homomorphism $N\to\R$
whose tangent homomorphism annihilates $H$, $\R^+=[0,\infty)$.
We shall prove that
$\gR(C)=N^+_\chi$ for some $\chi$ if the degree of contact of $C$ with
certain subspaces is sufficiently high.
The role of the degree of contact was notified in \cite{Gi} where global
$\Ad$-invariant cones were characterized --  while the final answer was
formulated by the algebraic language, in fact,
the globality of an invariant cone in a Lie algebra is determined
by the degree of contact of the cone with the linear sum of two distinct
nilpotent subalgebras.
In this article we use Carnot--Caratheodory metrics to prove the result
mentioned above. Probably, these metrics can be a natural and essential
tool in Geometric Control Theory, in particular, for the investigation
of attainable  sets. In any way, the usage of Carnot--Caratheodory metrics
clarifies the dependence of $\gR(C)$ on the degree of contact. They also
give quantitative versions for the criterion of controllability (\ref{ccrit}).
For a discussion of the role of Lie groups and algebras in Control Theory
and further references, see \cite{Br}, \cite{CKM}.
\section{Preliminaries and statement of the result}
\ssun{Realization of simply connected nilpotent Lie groups}
Let $\cN$ be a nilpotent Lie algebra.
The corresponding Lie group $N$ can be realized as $\cN$ with
the group multiplication defined by the Campbell--Hausdorff formula:
\begin{equation}\label{CHF}
xy=x+y+{\frac{1}{2}}[x,y]+P_3(x,y)+\dots +P_d(x,y),
\end{equation}
where $P_k(x,y)$, $k=3,\dots,d$, is the sum of Lie products of the
length $k$.
Thus $P_k$ is a homogeneous polynomial of the degree $k$ with values in $\cN$.
It follows from (\ref{CHF}) that for all $\xi\in\cN$
\begin{equation}\label{exp}
\exp(\xi)=\xi,\ins{}\inv\xi=-\xi,
\end{equation}
$e=0$, and the multiplicative commutator has the form
\begin{equation}\label{comm}
\{x,y\}= xyx^{-1}y^{-1}=[x,y]+\ (\hbox{Lie products of the length}>2)
\end{equation}
We shall consider simultaneously the Lie group
and the Lie algebra structures. In particular, we keep
the vector notation for addition
and multiplication by scalars. It will be convenient to fix the euclidean
distance in $\cN$ which will be denoted by $\abs\ $. Thus, for example,
\begin{equation}\label{mua}
x^n=nx\ins{and}\abs{x^n}=\abs n\abs x\ins{for all integer}n.
\end{equation}
Put
\begin{equation}\label{cent}
\cN^1=\cN,\quad\cN^{k+1}=[\cN,\cN^k];\quad\quad
\cN^1\supset\cN^2\supset\dots\supset\cN^d\supset\cN^{d+1}=\{0\},
\end{equation}
where $\cN^d\neq\{0\}$. Then $\cN^{k}$, $k=1,\dots,d$ is also a normal
subgroup of $\cN$ which sometimes will be denoted by $N^k$.
For each $k$ chose in $\cN^k$ a complementary to $\cN^{k+1}$ subspace
$\cN_k$. Then $\cN$ is the linear direct sum of these subspaces
\begin{equation}\label{grad}
\cN=\cN_1\oplus\cN_2\oplus\cdots\oplus\cN_d
\end{equation}
and $\cN_1$ generates $\cN$ as a Lie algebra.
\ssun{Graded nilpotent Lie algebras} If
$[\cN_k,\cN_l]\subseteq\cN_{k+l}$
then (\ref{grad}) is the gradation of $\cN$ which satisfies the
additional condition
\begin{equation}\label{equal}
[\cN_1,\cN_l]=\cN_{l+1},\quad l=1,\dots,d.
\end{equation}
which is equivalent to the assumption that $\cN_1$ generates $\cN$.
Let us pick $x\in\cN$, decompose it according to (\ref{grad})
$$x=x_1+\dots+x_d,$$
and put
$$Dx=x_1+2x_2+\dots+dx_d.$$
If (\ref{grad}) is the gradation then $D$ is a differentiation of $\cN$ and
\begin{equation}\label{det}
\de_t:\,x_1+x_2+\dots+x_d\ \to \ tx_1+t^2x_2+\dots+t^dx_d,\quad t>0
\end{equation}
is the corresponding one-parametrical group of automorphisms written
in the multiplicative form. Since $\exp$ is identical, $\de_t$ is
also the isomorphism of the group $N$. Further,
the group $\{\de_t\}_{t\geq0}$ can be extended to the complexification of
$\cN$ and nonzero complex values of $t$. This implies that
\begin{equation}\label{dem}
\de_{-1}:\,x_1+x_2+\dots+x_d\ \to \ -x_1+x_2+\dots+(-1)^kx_k+
\dots+(-1)^dx_d
\end{equation}
is an isomorphism of $\cN$ and $N$.
\ssun{Asymptotic group} If (\ref{grad}) is not a gradation then
the formula
\begin{equation}\label{lim}
[x,y]^a=\lim_{t\to\infty}\inv\de_t[\de_tx,\de_ty]
\end{equation}
defines the asymptotic Lie product in $\cN$ and (\ref{grad})
is the gradation for it.
This gradation also can be defined
by the standard factorization procedure for the filtration (\ref{cent}):
the new Lie bracket for $x\in\cN_k$, $y\in\cN_l$ is
the projection of the old one in $\cN^{k+l}$ to $\cN_{k+l}$ along
$\cN^{k+l+1}$. Indeed, let
$$[x,y]=[x,y]_1+\dots+[x,y]_d,\quad
x=x_1+\dots+x_d,\quad
y=y_1+\dots+y_d$$
be corresponding to (\ref{grad}) decompositions; then
\begin{equation}\label{dex}
\inv\de_t[\de_t x,\de_t y]=\sum_{p\geq k+l} t^{k+l-p}[x_k,y_l]_p=
\sum_{p=k+l}[x_k,y_l]_p+\al(t),
\end{equation}
where $|\al(t)|=O(\fracm1t)$,
and the limit can be easily calculated in this notation.
The corresponding group $N_a$ also can be realized as a limit.
For each $t>0$, put
$$x\mult t y=\de_t^{-1}(\de_tx\cdot\de_ty).$$
This introduces in $\cN$ the structure of a Lie group isomorphic to $N$.
It follows from the Campbell-Hausdorff formula and (\ref{dex}) that
for every $x,y\in N$ there exists the limit
\begin{equation}\label{grm}
x\mult a y=\lim_{\sssize t\to\infty}x\mult t y
\end{equation}
Moreover, by (\ref{dex}) and (\ref{CHF})
\begin{equation}\label{rest}
x\mult t y=\de_t^{-1}(\de_tx\cdot\de_ty)=x\mult a y+\be(x,y,t),
\end{equation}
where
$$|\be(x,y,t)|\leq \frac{A}t\,(|x|+|x|^d)(|y|+|y|^d)$$
and $A>0$ depends only on the algebra $\cN$.
\ssun{A construction for curves in $\gT(C)$} Let $\ga_1:\,[0,a_1]\to N$ and
$\ga_2:\,[0,a_2]\to N$ be paths in $N$ starting at $e$: $\ga_k(0)=e$,
$k=1,2$.
Put
$$\ga_1\cdot\ga_2(t)=\left\{\begin{array}{cc}\ga_1(t),&t\in[0,a_1]\cr
\ga_1(a_1)\ga_2(t-a_1),&t\in[a_1,a_1+a_2]\end{array}\right.$$
Then, for any $C\subseteq\cN$,  $\ga_1\cdot\ga_2$ belongs to $\gT(C)$
if so are $\ga_1$ and $\ga_2$. An important particular case of this
construction is as follows: for $\xi\in\cN$ put $\ov\xi(t)=\exp(t\xi)$,
$t\in[0,1]$, and
\begin{equation}\label{cur}
\ov x=\ov\xi_1\cdot\ov\xi_2\cdot\dots\cdot\ov\xi_n,\ins{where}
x=(\xi_1,\dots,\xi_n),\ \xi_1,\dots,\xi_n\in C.
\end{equation}
Then $\ov x(n)=\exp(t\xi_1)\dots\exp(t\xi_n)$. Furthermore, if
$\xi_1,\dots,\xi_n\in C$ then $\ov x\in\gT(C)$.
For the Riemannian left invariant metric defined by the euclidean norm
$\abs\ $
\begin{equation}\label{add}
\La(\ov x)=\La(\ov\xi_1)+\dots+\La(\ov\xi_n),
\end{equation}
where $\La(\ga)$ denotes the length of the curve $\ga$.
The length of the curve $\ov\xi$, $\xi\in C$, can be easily derived:
$\La(\ov\xi)=|\xi|$. Hence for $\xi_1,\dots,\xi_n\in C$
\begin{equation}\label{lxx}
\La(\ov x)=\abs{\xi_1}+\dots+\abs{\xi_n}.
\end{equation}
Clearly, each curve in $\gT(C)$ can be approximated by curves of the type
(\ref{cur}).
\ssun{\ccm{}s}
Recall the definition of \ccm{}s which also
are known as subriemannian or nonholonomic Riemannian ones. Any
euclidean norm $\abs\ \,$ on $\cN_1$ uniquely
determines the left invariant norm on the left invariant distribution of
subspaces generated by $\cN_1$.
Hence the length of a curve $\ga\in\gT(\cN_1)$ can be defined
by the standard formula
\begin{equation}\label{length}
\La(\ga)=\int_a^b\abs{\ga'(t)}_{\ga(t)}\,dt.
\end{equation}
Since $\ov\xi\in\gT(\cN_1)$ for $\xi\in\cN_1$,
the left invariance of the \ccm{} implies that (\ref{add}) and
(\ref{lxx}) are
true for curves described in Subsection~2.4 with $C=\cN_1$.
The Carnot--Caratheodory distance $\ka(x,y)$ between $x$ and $y$ is
defined as the least lower bound for lengths of curves in  $\gT(\cN_1)$
which join $x$ and $y$:
\begin{equation}\label{dist}
\ka(x,y)=\inf\{\La(\ga):\,\ga\in\gT(\cN_1),\ga:[a,b]\to N,\ga(a)=x,\ga(b)=y\}.
\end{equation}
Each piecewise smooth curve in $\cN_1$ has
the unique lift to the curve in $\gT(\cN_1)$. Hence the natural projection
$\pi_k:\,N\to N/N^k$, $k=2,\dots,d$, keeps the length of a curve $\ga$
in $\gT(\cN^1)$ which is equal to  the usual euclidean length
of $\pi_2\ga$:
\begin{equation}\label{eqlen}
\La(\pi_2\ga)=\La(\pi_3\ga)=\dots=\La(\pi_d\ga)=\La(\ga).
\end{equation}
If (\ref{grad}) is the
gradation satisfying (\ref{equal}) then $\de_t$ is an automorphism,
hence it commutes with the lifting procedure. Therefore, $\de_t$
in (\ref{det}) is a {\it metric dilation} for $t>0$, i.e.
\begin{equation}\label{dil}
\rho(\de_t(x),\de_t(y))=t\rho(x,y),\ins{}x,y\in N.
\end{equation}
Since $\de_t(\xi)=t\xi$ for $\xi\in\cN_1$, (\ref{dil})  follows from
(\ref{length}) and the definition of $\ka$. Further, (\ref{length}) and
(\ref{dem}) implies that $\de_{-1}$ is an isometry:
\begin{equation}\label{iso}
\ka(\de_{-1}(x),\de_{-1}(y))=\ka(x,y)\ins{for all}x,y\in N.
\end{equation}
If $\cN$ is not graded then $\de_t$ is not an automorphism but
this is true for the limit group defined by (\ref{grm}).
It can be equipped with the limit metric
\begin{equation}\label{mas}
\ka_a(x,y)=\lim_{t\to\infty}\ka_t(x,y),
\ins{where}\ka_t(x,y)=\fracm1t\ka(\de_tx,\de_ty).
\end{equation}
The asymptotic group could be realized as
the Gromov--Hausdorff limit of metric spaces -- groups with left
invariant metrics $\ka_t$. For more details on this subject, see
(\cite{Gr}), (\cite{VeGe}).
\ssun{Inner metrics}
We use a definition of the inner metric which
is  equivalent to the standard one in the class of left invariant metrics
on Lie groups (see \cite{BPS}). Let $\rho$ be a left invariant metric which
is compatible with the topology, $B(x,r)$ denote the open ball at $x$ of
radius $r>0$, $B(x,0)=\{x\}$. Put $B(r)=B(e,r)$ and $B(0)=\{e\}$. The left
invariance of the metric $\rho$ means that
$$B(x,r)=xB(r).$$
We shall say that $\rho$ is
{\it inner} if
\begin{equation}\label{inn}
B(r)B(s)=B(r+s)\ins{for all}r,s\geq0.
\end{equation}
The product of sets is taken pointwise: $AB=\{ab:\,a\in A,\,b\in B\}$.
The same equality is true for closed ball since they are compact. Taken
together with the left invariance, this implies for any two points in $N$
the existence of the shortest curve which joins them.
If $H$ is a normal closed subgroup of a Lie group $G$ with the inner metric
$\rho$ and $\pi:\,G\to G/H$ is the canonical projection then
$\td B(r)=\pi B(r)$ is the unit ball for the metric
\begin{equation}\label{sbm}
\td\rho(\td x, \td y)=\inf\{\rho(xh_1,yh_2):\,h_1,h_2\in H\}.
\end{equation}
Clearly, this metric is inner because (\ref{inn}) is satisfied for it.
Furthermore, Riemannian and \ccm{}s are inner -- this follows
from (\ref{dist}) for \ccm{}s and from the analogous formula for Riemannian
ones.
The identity (\ref{inn}), in particular, implies that all inner
metrics are equivalent "in large": for each pair of inner metrics
$\rho$, $\rho'$ which define the same topology and any $\ep>0$, there exist
$C,c>0$ such that
\begin{equation}\label{equ}
x,y\in G,\ \rho(x,y)>\ep\Longrightarrow
c<{{\rho'(x,y)}\over{\rho(x,y)}}<C
\end{equation}
Indeed, for some $C,c>0$ inclusions
$B(c\ep)\subseteq B'(\ep)\subseteq B(C\ep)$
holds, where $B'(\ep)$ is the $\rho'$-ball. By (\ref{inn}),
$$B(nc\ep)\subseteq B'(n\ep)\subseteq B(nC\ep)$$
for all positive integer $n$, and the left invariance of these metrics
implies (\ref{equ}).
Thus, each inner metric asymptotically (for great distances) equivalent to a
\ccm{}. Note that the \ccm{} for the asimptotic group is self-similar -- it
admits metric dilations. Left invariant inner metrics on topological groups
were studied in (\cite{Be1}), (\cite{Be2}). For Lie groups, they are
Finsler (maybe nonholonomic) ones.
\ssun{Degree of contact}
Let  $\eta $ be an increasing function defined on some
interval in  $\R$ with the left endpoint  $0$,
 $\lim_{\ep\to 0}\limits\eta(\ep)=0$, $\cL$ be a linear subspace of
the euclidean space $\cN$. We shall say that a cone  $C$ has the
{\it degree of contact with $\cL$ at the point  $x\in C$  greater than}
$\eta$ if
\begin{equation}\label{deg}
\dist(C,x+y)=o(\eta(\abs y))\ins{as}y\to 0\ins{in}\cL
\end{equation}
The  degree of contact of $C$ with $\cL$  at  $x$ is
{\it greater or equal to}  $\eta$ if there exist $Q>0$
and a neighborhood  $U$ of zero in  $\cL$ such that
$$Q\eta(\abs y)\geq\dist(C,x+y)\ins{for all}y\in U$$
Suppose that  $x\notin\cL$; then the degree of contact is {\it equal}
to $\eta$   if it is greater or equal and the inverse inequality
holds with some another constant.
If  $x\in\cL$ then one has to replace  $\cL$ in this definition
to any subspace  $\cL'\subset\cL$ complementary to  $\R x$ in $\cL$,
the definition doesn't depend on the choice of  $\cL'$.

If  $\eta(\ep)=\ep^a$  then  $a$ will be called the degree of contact
and denoted by   $\cont(C,\cL,x)$;  $\cont(C,\cL,x)>a$
($\cont(C,\cL,x)\geq a$, $\cont(C,\cL,x)=a$)  will mean that the degree
of contact is greater than (respectively, greater or equal to, equal to)
$\ep^a$. For example, the degree of contact of any Lorentian cone with
each its tangent hyperplane is equal to 2: $\cont(C,T_x\partial C,x)=2$
for any $x\in\partial C$.
\ssun{Statement of the main result}
The result is proved for slightly more general setting then it was
mentioned above: the group need not be nilpotent in general -- it is supposed
to be a semidirect product of $\R$ and a nilpotent group.
We keep the notation for $\cN$, $N$,
particularly (\ref{cent}), which were introduced above.
Let $\R^+$ denote the set of all nonnegative real numbers. By
$\cG$ we denote a real Lie algebra of a simply connected Lie group
$G$ which is a semidirect product of $\R$ and a nilpotent group $N$,
$\chi:\,G\to\R$ be the projection homomorphisms to the factor $\R$.
Set $G^+=\inv\chi(R^+)$, and $\cG^+=\inv d_e\chi(\R^+)$ (hence
$\exp(\cG^+)=G^+$).
\begin{theorem} Let $C\subseteq\cG^+$ be a convex closed generating cone,
$d>1$, $p\in\partial C\cap\cN^d$; further,
suppose that there exists $v\in\cZ(p)\cap\Int(C)\neq\empty$ admitting
$\ad(v)$-invariant linear subspace $\cN_1\subset\cN^1$ complementary
to $\cN^2$, and
\begin{equation}\label{cont}
\cont(C,\cN_1,p)>\frac{d}{d-1}.
\end{equation}
Then $\gR(C)=G^+$.
\end{theorem}
\section{Quantitative versions of the controllability}
In the following two lemmas the Lie algebra $\cN$ is supposed to be graded
as in (\ref{grad}), (\ref{equal}).
\begin{lemma}
Let $\rho$ be an inner metric in $N$, $z\in N_d$, $\ep>0$, and $r=\rho(e,z)$.
Then
\begin{equation}\label{qccrit}
e\in B(z,\ep)^n\ins{\iff}n>\left(\frac r\ep\right)^{{d}\over{d-1}}.
\end{equation}
Moreover, for any $s>0$, if
\begin{equation}\label{ttt}
n^{\frac1d}r+s<n\ep
\end{equation}
then $B(z,\ep)^n\supset B(s)$.
\end{lemma}
\proof Since $z$ belongs to the center of $N$,
\begin{equation}\label{ball}
B(z,\ep)=zB(\ep)=B(\ep)z.
\end{equation}
Taken together with (\ref{inn}) and (\ref{mua}), \lasteqn0 implies that
\begin{equation}\label{power}
B(z,\ep)^n=(zB(\ep))^n=(nz)B(n\ep)=B(nz,n\ep).
\end{equation}
Therefore, the inclusion in the left side of (\ref{qccrit})  is equivalent
to the inequality
\begin{equation}
\rho(e,nz)<n\ep.
\end{equation}
By (\ref{dil}) and (\ref{det}),
$$\rho(e,nz)=n^{\frac1d}\rho(e,z)=n^{\frac1d}r.$$
Thus inequality \lasteqn1 holds \iff{}
$$n^{\frac1d}r<n\ep.$$
This is equivalent to the right part of (\ref{qccrit}).

For $w\in N$, the assumption
$\rho(w,nz)<n\ep$
is equivalent to $w\in B(nz,n\ep)$. By the triangle
inequality,
$$\rho(e,nz)<t,\ \rho(e,w)<s\ \Longrightarrow\  w\in B(nz,t+s),$$
hence the inclusion $B(nz,t+s)\supseteq B(s)$.
Put $t=n^{\frac1d}r$. Then, according to (\ref{power}), \lasteqn0,
and (\ref{ttt}), we receive the desired inclusion.\qed
\begin{corollary} For non-graded $\cN$, if $s>0$ then there exists $C>0$ such
that for all $\ep>0$
\begin{equation}
C(n^{\frac1d}r+s)<n\ep
\end{equation}
implies $B(z,\ep)^n\supset B(s)$.
\end{corollary}
\proof This follows from the existence of the asymptotic metric (\ref{mas})
for which the assertion holds by the lemma, and (\ref{equ}).\qed
\begin{lemma}\label{d-1}
Let $x\in\cN_{d-1}$, $\ep>0$. Suppose that $B(x,\ep)^n\cap\cN_d\neq\{0\}$.
Then $e\in B(x,2\ep)^{2n}$.
\end{lemma}
\proof Since the metric is left invariant, $B(x,\ep)=xB(\ep)$.
Let
$$z=xy_1xy_2\dots xy_n\in\cN^d$$
where $y_1,y_2\dots y_n\in B(\ep)$.
If $d$ is odd then $\de_{-1}x=x$, and $\de_{-1}z=-z=\inv z$. Hence
$$e=z\inv z=
xy_1\dots xy_nx\de_{-1}(y_1)xy_2\dots x\de_{-1}(y_n)\in B(x,\ep)^{2n}.$$
If $d$ is even then $\de_{-1}x=\inv x=-x$, and $\de_{-1}z=z$. Therefore,
$$e=z\de_{-1}(\inv z)=
xy_1\dots xy_n\de_{-1}(y_n)x\dots\de_{-1}(y_1)xe
\in B(x,2\ep)^{2n}$$
since $y_n\de_{-1}(y_n)\in B(2\ep)$.\qed
\begin{corollary}
Let $x\in\cN^{d-1}$, $\ep>0$. Then
\begin{equation}
n>2\left(\frac{2r}{\ep}\right)^{\frac d{d-1}}
\quad\Longrightarrow\quad e\in B(x,\ep)^{n}.
\end{equation}
\end{corollary}
\proof Applying (\ref{qccrit})  to the factor group $N/N^d$ we receive
$$n>\left(\frac{r}{\ep}\right)^{\frac d{d-1}}
\quad\Longrightarrow\quad B(x,\ep)^{n}\cap N^d\neq\empty,$$
and Lemma~\lastlem1, with $\ep$ replaced by $\frac\ep2$, implies
the desired inclusion.\qed
There is a natural way to realize
any finite dimensional nilpotent Lie algebra $\cN$ as
a factor algebra of a finite dimensional graded Lie algebra $\wt\cN$.
Let $x_1,\dots,x_l$ be a set of generators for $\cN$ (the linear
basis of $\cN_1$). Then $\cN$ is the homomorphic image
of the free Lie algebra $\cF$ generated by $x_1,\dots,x_l$. The kernel of the
homomorphism includes the ideal $\cF^d$ generated by all products of
length $>d$. This means that $\cN$ is the homomorphic image of the finite
dimensional Lie algebra $\wt\cN=\cF/\cF^d$ whose natural gradation
satisfies the condition (\ref{equal}). Let $\pi:\,\wt\cN\to\cN$
denote this homomorphism. Clearly, $\pi\wt\cN^k=\cN^k$ for $k=1,\dots,d$.
Note that $\wt\cN$ has the same height $d$ and that generating
spaces $\cN_1$ and $\wt\cN_1$ may be identified.
Thus the euclidean norm
$\abs\ $ in $\cN_1$ defines \ccm{}s $\ka$ and $\td\ka$ in $N=\cN$ and
$\wt N=\wt\cN$ respectively.
We shall equip with \ $\td{}$ \ symbols denoting objects in $\wt\cN$
corresponding to objects in $\cN$. Put
$$\ka(x)=\ka(e,x),\quad\td\ka(\td x)=\td\ka(\td e,\td x).$$
Clearly, $\ka(\pi\td x)\leq\td\ka(\td x),\quad\td x\in\wt N$.
In the following theorem we do not assume that $\cN$ is graded but keep
the notation of the previous section.
\begin{theorem} Let $d>2$, $k=d$ or $k=d-1$, $x\in N^k\setminus N^{k+1}$,
and $r=\ka(x)$. Then there exists $Q\geq1$ such that for all $\ep>0$
the condition
\begin{equation}\label{nnn}
n>Q\left(\frac r\ep\right)^{\frac{k}{k-1}}
\end{equation}
implies that $e\in B(x,\ep)^n$, where $B(x,\ep)$ is the ball for the
\ccm{} $\ka$.
\end{theorem}
\proof
If $\cN$ is graded as in (\ref{grad}), (\ref{equal}), then the assertion
of the theorem is an easy
consequence of Lemma~1 for $k=d$ and Corollary~2 for $k=d-1$. In general
case, let us realize $\cN$ as the factor algebra of the graded Lie
algebra $\wt\cN$ by the construction described above. Let
$\pi\td x\in\td\cN^k\setminus\cN^{k+1}$, $\pi\td x=x$, and
put  $\td r=\td\ka(\td x)$, $K=\frac{\td r}{r}$. Then, by  Lemma~1
or Corollary~2, there exists $A>0$ such that
$$n>A\left(\frac{\td r}\ep\right)^{\frac{k}{k-1}}\ins{implies}
\td e\in\td B(\td x,\ep)^n.$$
Since $\pi\td B(\td x,\ep)^n=B(x,\ep)^n$,  the inclusion  $e\in B(x,\ep)^n$
is true for 
$$Q=\max\{1,K^{\frac{k}{k-1}}A\}.\qedm$$
The following theorem is in fact a reformulation of Theorem~2 by another
words. Let $\rho$ be the Riemannian metric defined by the euclidean norm
$|\ |$ in $\cN$ and $\ka$ be the \ccm{} corresponding to the restriction
of this norm to $\cN_1$. Put $\cB(r)=\{\xi\in\cN:\,|\xi|<r\}$
\begin{theorem} Let $d$, $k$, $x$, $\rho$ be as above, and let $r=\rho(e,x)$.
Then there exists $P>0$ such that for all safficiently small $\ep>0$
the group $N$ admits a closed curve $\ga\in\gT(x+\cB(\ep))$
whose $\rho$-length satisfies the inequality
\begin{equation}\label{leq}
\La_\rho(\ga)\leq P\left(\frac r\ep\right)^{\frac{k}{k-1}}.
\end{equation}
\end{theorem}
\proof  Put $R=|x|$; clearly, $r\leq R\leq\ka(x)$.
If $\ep>R$ then the assertion is evident. Hence we may assume that $\ep<R$.
Then there exists $a\in(0,1)$
such that $x+\cB(\ep)$ includes the $\rho$-ball at $x$ of radius $a\ep$
for all $\ep\in(0,R)$ (recall that we identify $\cN$ and $N$). Hence
it includes the $\ka$-ball $B(x,a\ep)$.
Let $Q$ be as in Theorem~2. Then there exists
integer $n$ which satisfies inequalities
$$Q\left(\frac{\ka(x)}{a\ep}\right)^{\frac{k}{k-1}}<n\leq
2Q\left(\frac{\ka(x)}{a\ep}\right)^{\frac{k}{k-1}}.$$
Then, by the first of them and Theorem~2, there exist
$x_1,\dots,x_n\in B(x,a\ep)$ such that  $x_1\dots x_n=e$.
Since $\exp$ is identical, it follows from the construction of
Subsection~2.4 that $\ov x_1,\dots\ov x_n\in\gT(B(x,a\ep))$ and the curve
$\ga=\ov x_1\cdot\dots\cdot\ov x_n$ is closed.
By (\ref{lxx}),
$$\La_\rho(\ga)=|x_1|+\dots+|x_n|\leq2Rn<
P\left(\frac r\ep\right)^{\frac{k}{k-1}},$$
where
$$P=
4QR\left(\frac{\ka(x)}{ar}\right)^{\frac{k}{k-1}}.$$
This proves the theorem.\qed
\section{Attainable  sets} Everywhere in this section we suppose that
the assumption of Theorem~1 are satisfied.
Let $\cG$ be equipped with
the euclidean norm $\abs\ $ and $G$ with the corresponding left invariant
Riemannian metric $\rho$ and $N$ with \ccm{} $\ka$. Then the semidirect
product $G=\R\smdr N$ is defined by the one-parametrical
group $A_t=e^{t\ad(v)}$, $t\in\R$, of group automorphisms of $N$.
Since $\exp$ is identical for the coordinate system in $N$ which we use,
$A_t$ is also the one-parametrical group of automorphisms of $\cN$. Hence
$A_t$ is linear in these coordinates. Put
\begin{equation}\label{M}
M=\sup\{\norm{A_t}:\,\abs t\leq1\}.
\end{equation}
The multiplication law in the group $G$ can be written explicitly:
$$(t,x)(s,y)=(t+s,(A_{-s}x)y),\ins{where}x,y\in N,\ t,s\in\R.$$
We denote $v=(1,0)$.
By the assumtion of the theorem,
\begin{eqnarray}
p+v\in\Int(C),\label{p+v}\\
v\in\cZ(p).\label{v}
\end{eqnarray}
Set
$$\cB=\{\xi\in\cG:\,\abs\xi<1\},\quad\cB_1=\cB\cap\cN_1,$$
and let $B_\ka(r)$ the $\ka$-ball with the center $e$ of the radius $r$.
In the following lemma we consider these sets as subsets of the group $N$.
\begin{lemma} For any $r>0$
\begin{equation}
B_\ka(r)=\Cup_{n=1}^{\infty}\left(\fracm rn\cB_1\right)^n.
\end{equation}
\end{lemma}
\proof For each $x\in\cB_1$ the curve $\ov x$ belongs to $\gT(\cB_1)$.
Hence $\fracm rn\cB_1\subseteq B_\ka(\fracm rn)$. By (\ref{inn}),
the left side of the equality includes the right one. Let
$\ga:\,[0,r]\to N$ be a curve in $\gT(\cB_1)$.
It follows from the definition of \ccm{} that the open ball $B_\ka(r)$
is filled by points $\ga(r)$ for such curves $\ga$. Let $\la_g(h)=gh$ be
the left shift by $g$. The endpoint of the curve
$$\ga_n=\ov x_1\cdot\dots\cdot\ov x_n,$$
where
$$
x_k=\fracm1n d_{\ga(t_k)}\inv\la_{\ga(t_k)}(\ga'(t_k)),
\quad t_k=\fracm{rk}{n},\quad k=1,\dots,n,$$
belongs to the right side of \lasteqn0. Clearly, $\ga_n(t)\to\ga(t)$ as
$n\to\infty$ for each $t\in[0,r]$. Hence
the right side of \lasteqn0 is dense in the left one. Further, it follows
from (\ref{comm}) and (\ref{exp}) that $(s\cB_1)^k$ is open for
sufficiently large $k$ depending only on $\cN$. Therefore, the right side of
\lasteqn0 is open; let us denote it by $\td B(r)$. By standard arguments
it is not difficult to show that
$$\Cup_{n=1}^{\infty}\left(\fracm rn\cB_1\right)^n=
\Cup_{n=1}^{\infty}\left(\fracm r{2^n}\cB_1\right)^{2^n}.$$
Hence $\td B(\fracm r2)^2=\td B(r)$, and this division procedure can
be continued. This implies that $\td B(r)$ coincides with the interior
of it's closure. Since $\td B(r)$ is dense in $B_\ka(r)$ and
open, $\td B(r)=B_\ka(r)$.\qed
\begin{lemma}
If $q\in N^d$, $\ep>0$, $t>0$, $n\in\N$, and $nt<1$ then
\begin{eqnarray}
(t,q+M\ep\cB_1)^n\supseteq (nt,nq+(\ep\cB_1)^n);\label{cka}\\
(t,B_\ka(q,M\ep))^n\supseteq(nt, B_\ka(q,\ep)^n)\label{pka}.
\end{eqnarray}
\end{lemma}
\proof Let $x_1,\dots,x_n\in(p+\ep\cB_1)$.
Put $\td x_k=(t,A_{(n-k)t}x_k)$. Then
\begin{eqnarray*}
\td x_1\td x_2\dots\td x_n=(t,A_{(n-1)t}x_1)(t,A_{(n-2)t}x_2)\dots(t,x_n)=\\
(2t,A_{(n-2)t}(x_1x_2))\dots(t,x_n)=\dots=(nt,x_1x_2\dots x_n).
\end{eqnarray*}
Since $\cN_1$ is $A_t$-invariant, $nt<1$, by (\ref{v}) and (\ref{M}),
$A_{(n-k)t}x_k\in(p+M\ep\cB_1)$ for all $k=1,\dots,n$. To prove (\ref{cka}),
it remains to note that $(q+\ep\cB_1)^n=nq+(\ep\cB_1)^n$ because $p$ belongs
to the center of $N$. The inclusion (\ref{pka}) follows from the same
equality, with $x_k\in B_\ka(q,\ep)$, and the inequality
$\ka(e,A_tx)\leq M\ka(e,A_tx)$ which is an easy consequence of the definition
\ccm{} $\ka$.\qed
The following elementary lemma whose assertion could be a definition
of the degree of contact was already proved in (\cite{Gi}).
We omit the proof -- it is rather long than hard.
Let $\cL$ be as in Subsection~2.7. Put
$$
\cB_\cL=\{y:\,y\in\cL,\,\abs y\leq1\}.
$$
\begin{lemma}
Let $C$ be a generating closed cone in $\cN$, $x\in C$, $x\neq0$.
Suppose that  $\cont(C,\cL,x)>a\geq1$, $v\in\cE$, and $x+v\in\Int (C)$.
Then there exists a function $\vf$ defined on some
interval $(0,\al)$, $\al>0$, such that $\vf(\ep)=o(\ep^a)$ as $\ep\to0$ and
$$x+\vf(\ep)v+\ep\cB_\cL\subset\Int(C)$$
for all $\ep\in(0,\al)$.\qed
\end{lemma}
\par\noindent{\it Proof of Theorem 1.\kern3pt} Let $B_\ka(p,\ep)$ be the
Carnot--Caratheodory ball in $N$, $B_\ka=B(e,1)$, $\al$, $\vf$ be as in
Lemma~5 with $a=\frac{d}{d-1}$, $\cL=\cN_1$. There exists a function
$\psi$ such that
\begin{eqnarray}
\lim_{t\to0}\frac{\vf(t)}{\psi(t)}=0,\label{aaa}\\
\ins{}\psi(t)=o(t^{\frac{d}{d-1}})\label{bbb}.
\end{eqnarray}
It follows from Corollary~1 that there exists $A>0$ such that
$$n>A\ep^{-\frac{d}{d-1}}\quad\Rightarrow\quad B_\ka(p,\ep)^n\supseteq B_\ka.
$$
For these $n$ and sufficiently small $\ep>0$, applying (\ref{pka}) we receive
$$(\psi(\ep), B_\ka(p,M\ep))^n\supseteq(n\psi(\ep), B_\ka(p,\ep)^n)
\supseteq(n\psi(\ep), B_\ka).$$
By (\ref{bbb}), since $n$ can be chosen satisfying the inequality
$n<2A\ep^{-\frac{d}{d-1}}$,
the set $\gR(C)$ includes the ball $B_\ka$, hence the group $N$ and the
halfspace $G^+$. Therefore, it is sufficient to prove the inclusion
\begin{equation}
\gR(C)\supset (\psi(\ep), B_\ka(p,M\ep))
\end{equation}
for $\ep\in(0,\al)$ for some $\al>0$. It follows from Lemma~5 and
(\ref{aaa}) that
$$C\supset (\psi(\ep)v,p+M^2\ep\cB_1)$$
if $\ep$ is sufficiently small. Since $C$ is a cone,
$$C\supset\fracm1n(\psi(\ep)v,p+M^2\ep\cB_1),\quad n\in\N,$$
hence
$$\gR(C)\supset\left(\fracm1n(\psi(\ep), p+M^2\ep\cB_1)\right)^n\supseteq
\left((\psi(\ep),p+(\fracm{M\ep}{n}\cB_1)^n\right)$$
by (\ref{cka}), and the desired inclusion \lasteqn0 follows from
Lemma~3.\qed

\bibliographystyle{unsrt}

\vspace{5mm}
%
%
\textsc{Omsk State University, prosp. Mira 55a, Omsk 644077, Russia}

\textit{Current address}:
\textrm{Mathematical Department, Omsk State University,
prosp. Mira, Omsk, Russia}

\textit{E-mail address}:
\texttt{gichev@math.omsu.omskreg.ru}

\vspace{5mm}
%
%
\end{document}